\documentclass[12pt]{amsart}

\usepackage{amssymb}
\usepackage{amsmath}
\usepackage{latexsym}
\usepackage{epsfig}
\usepackage{graphics}
\usepackage{color}

\newcommand{\A}{\ensuremath{\mathbb{A}}}
\newcommand{\Q}{\ensuremath{\mathbb{Q}}}
\newcommand{\Z}{\ensuremath{\mathbb{Z}}}
\newcommand{\C}{\ensuremath{\mathbb{C}}}
\newcommand{\R}{\ensuremath{\mathbb{R}}}

\newcommand{\G}{\ensuremath{\mathbb{G}}}

\DeclareMathOperator{\Spec}{Spec}

\DeclareMathOperator{\divisor}{div}

\newcommand\cO{\mathcal{O}}

\newtheorem{lemma}[subsection]{Lemma}
\newtheorem{definition}[subsection]{Definition}  

\newtheorem{remark}[subsection]{Remark}
\newtheorem{proposition}[subsection]{Proposition}

\newtheorem{namedTheorem}[subsection]{\theoremname}
\newcommand{\theoremname}{testing}

\begin{document}

\title{SOME REMARKS ON TOROIDAL MORPHISMS}

\author[Denef]{Jan Denef}

\address[Denef]{University of Leuven, Department of Mathematics, Celestijnen-laan 200 B,
B-3001 Leuven (Heverlee), Belgium.}
\email{jan.denef@wis.kuleuven.be}


\date{\today}


\maketitle

\section{Introduction}

\noindent
This note contains some results related to the definition of toroidal morphisms over a field $k$ of characteristic zero. In \cite{Abramovich-Denef-Karu} this notion was defined by requiring that the base change of the morphism to an algebraic closure of $k$ is toroidal.
The notion of a toroidal morphism $f$ over an algebraically closed field was introduced long before by several authors, see e.g. \cite{Abramovich-Karu} and \cite{Cutkosky-Tor3folds}. Roughly the definition requires that for each closed point x of the source of $f$ one can choose formal toric coordinates at $x$ and formal toric coordinates at $f(x)$, such that in these coordinates the morphism is given by monomials. When $k$ is not algebraically closed, there is the natural question whether this remains true over the base field $k$ itself instead of over an algebraic closure of $k$. In this note we show that the answer is yes for toroidal morphisms between strict toroidal embeddings if the residue field $k(x)$ of $x$ equals $k$ or if $k(x)$ is algebraically closed. This is implied by Proposition \ref{prop2} below which is actually a stronger statement. An easy counterexample (Remark \ref{remark2-counterExample}) shows that the condition on $k(x)$ cannot be omitted. Proposition \ref{prop2} can be proved using Kato's paper \cite{Kato-LogStruct}, adapting the argument in section 3.13 of \cite{Kato-LogStruct}. However we preferred to provide a self-contained proof which does not use logarithmic geometry. Proposition \ref{prop2} in the special case of nonsingular toroidal embeddings is used in \cite{Denef-MonomializAndpAdicElimination} for applications of toroidalization to model theory. Proposition \ref{prop2} also implies that in the definition of toroidal morphisms, as formulated in \cite{Abramovich-Denef-Karu}, we can replace the completions by henselizations. This also holds for the definition of toroidal embeddings, see Remark \ref{remark0} below.

We will use without further mentioning the terminology and notation of \cite{Abramovich-Denef-Karu}, in particular we refer to \cite{Abramovich-Denef-Karu} for the notions of toroidal embeddings, strict toroidal embeddings, and toroidal morphisms. Moreover, $k$ will always denote a field of characteristic zero, except that the material in 1.1 up to \ref{lemma4} remains valid for any field $k$.

\section{Toroidal embeddings}
In \cite{Abramovich-Denef-Karu} an open embedding of algebraic varieties $U \subset X$ over $k$ is defined to be toroidal if its base change to an algebraic closure of $k$ is a toroidal embedding. Nevertheless we have the following proposition.

\begin{proposition}\label{prop1} Let $U \subset X$ be a strict toroidal embedding of varieties over $k$, and $x$ a closed point of $X$.Then there exist an affine toric variety $V$ over $k$ and an \'etale $k$-morphism $\varphi$ from an open neighborhood of $x$ in $X$ to $V$, such that locally at $x$ (for the Zariski topology) we have $U = \varphi^{-1}(T)$, where $T$ is the big torus of $V$.
\end{proposition}
\noindent This is proved in \cite{KKMS}, page 195, when $k$ is algebraically closed. However the proof remains valid in the general case due to Lemma \ref{lemma0} below.

\begin{definition} We call $(V,\varphi)$, with $V$ and $\varphi$ as in Proposition \ref{prop1}, an \emph{\'etale chart} for $U \subset X$ at $x$.
Note that there always exists an \'etale chart for $U \subset X$ at $x$ such that $\varphi(x)$ belongs to the closed orbit of $V$.
\end{definition}

\begin{lemma} \label{lemma0} Let $X$ be a normal algebraic variety over $k$ and $D$ a strict Weil divisor on $X$, i.e. the irreducible components of $D$ are normal. Let $x$ be a closed point of $X$ and $\widetilde\cO_{X,x}$ the completion of a strict henselization of $\cO_{X,x}$.
Let $y \in \widetilde\cO_{X,x}$ and assume that the Weil divisor $\divisor(y)$ of $y$ on $\Spec(\widetilde\cO_{X,x})$ is supported on the preimage $\widetilde{D}$ of $D$ in $\Spec(\widetilde\cO_{X,x})$. Then there exists $z \in \cO_{X,x}$ such that $y/z$ is a unit in $\widetilde\cO_{X,x}$.
\end{lemma}
\noindent \emph{Proof.} This is well known. Since we could not find a good reference, we include a proof. Because the natural morphism
$\Spec(\widetilde\cO_{X,x}) \rightarrow \Spec(\cO_{X,x})$ is flat, it induces a morphism $\tau$ from the group of Weil divisors on $\Spec(\cO_{X,x})$ to the group of Weil divisors on $\Spec(\widetilde\cO_{X,x})$, cf. \cite{EGA-4-4} Proposition 21.10.6. Moreover $\tau$ restricts to a bijection between the divisors supported on $D$ and those supported on $\widetilde{D}$, because the ideal in $\cO_{X,x}$ of any irreducible component of $D$ containing $x$ generates a prime ideal in $\widetilde\cO_{X,x}$, since these components are normal. Thus there exists a Weil divisor $W$ on $\Spec(\cO_{X,x})$ such that $\tau(W)=\divisor(y)$. The morphism $\tau$ induces an injection on the groups of divisor classes, indeed this follows easily by adapting the proof of Proposition 16 in section 1.10 of \cite{Bourbaki-CommAlg7}. Hence there exists $z \in \cO_{X,x}$ such that $\divisor(z)=W$. Thus $\tau(\divisor(z))=\divisor(y)$ and $y/z$ is a unit in $\widetilde\cO_{X,x}$. $\square$

\begin{remark}\label{remark0}\rm Proposition \ref{prop1} directly implies that in the definition of strict toroidal embeddings of varieties over $k$, as formulated in \cite{Abramovich-Denef-Karu}, we can replace the completions by henselizations. Hence this also holds for toroidal embeddings that are not necessarily strict. Indeed, for any toroidal embedding $U_X \subset X$ over $k$, and any closed point $x$ of $X$, there exists an \'etale morphism $f: X' \rightarrow X$ onto a neighborhood of $x$ such that $f^{-1}(U) \subset X'$ is a strict toroidal embedding. This assertion follows from Lemma \ref{lemma4} below with $A$ the strict henselization of $\cO_{X,x}$.
\end{remark}

\begin{lemma} \label{lemma4} Let $A$ be an excellent normal henselian local ring, and $P$ a prime ideal of height 1 in $A$. Let $\hat{A}$ be the completion of $A$. Assume for each height one prime ideal $P'$ of $\hat{A}$, with $P$ contained in $P'$, that $\hat{A}/P'$ is normal. Then $A/P$ is normal.
\end{lemma}
\noindent \emph{Proof.} Because $A/P$ is excellent, henselian and integral, its completion is also integral, by Corrolaire 18.9.2 in \cite{EGA-4-4} (or by Artin's Approximation Theorem \cite{Art} if A is moreover the strict henselization of the local ring of a closed point on a variety over $k$).
Thus $P\hat{A}$ is a prime ideal of $\hat{A}$. Since $P\hat{A}$ has height one (because $A$ and $\hat{A}$ are catenary), our assumption implies that $\hat{A}/P\hat{A}$ is normal. Thus the completion of $A/P$ is normal. By faithfully flat descent this implies that $A/P$ is normal (see e.g. Remark 2.24 in chapter 1 of \cite{Milne-EtaleCoh}). $\square$

\begin{definition} Let $U_X \subset X$ be a strict toroidal embedding of varieties over $k$. Denote by $D:= X \setminus U_X$ the toroidal divisor. The sheaf of \emph{logarithmic differential 1-forms} on $X$ is defined as the sheaf of $\cO_X$-modules $\Omega_X^1(\log D) := j_*(\Omega_{X_0}^1(\log D\cap X_0))$, where $j: X_0 \rightarrow X$ is any nonsingular open subscheme of $X$ with codimension $\geq 2$ such that $D \cap X_0$ is nonsingular.
\end{definition}
\noindent The sheaf $\Omega_X^1(\log D)$ of $\cO_X$-modules is locally free: as basis in a neighborhood of a closed point $x$ of $X$ one can take
$\frac{dx_1}{x_1}, \dots, \frac{dx_m}{x_m}$, where $x_1, \dots, x_m$ are the pullbacks, to the function field $K(X)$ of $X$, of the elements of a basis for the $\Z$-module of characters of the big torus of $V$, for any \'etale chart $(V, \varphi)$ for $U_X \subset X$ at $x$. Indeed, this follows from Proposition 15.5 in \cite{Danilov-ToricGeom} (the assumption there that $k=\C$ is not necessary).

\section{Logarithmically smooth morphisms}

\begin{definition} Let $U_X \subset X$ and $U_B \subset B$ be strict toroidal embeddings of varieties over $k$, and $x$ a closed point of $X$. Denote the toroidal divisors by $D_X:=X \setminus U_X$, $D_B:=B \setminus U_B$. Let $f: X \rightarrow B$ be a dominant $k$-morphism mapping $U_X$ into $U_B$. The morphism $f$ is called \emph{logarithmically smooth} at $x$ (with respect to $U_X \subset X$ and $U_B \subset B$) if the sheaf of $\cO_X$-modules
$$ \Omega_X^1(\log D_X) / f^*(\Omega_B^1(\log D_B))$$
is locally free at $x$. This is equivalent with the condition that the fiber at $x$ of this sheaf has dimension $\,\dim X -\dim B$ as vector space over the residue field $k(x)$ of $X$ at $x$. (Note that we required $f$ to be dominant.)
\end{definition}

\begin{remark}\label{remark1}\rm
Clearly, if $f: X \rightarrow B$ is toroidal with respect to $U_X \subset X$ and $U_B \subset B$, then $f$ is logarithmically smooth at each closed point of $X$. The converse is also true, this follows from Theorem 3.5 and Proposition 3.12 in \cite{Kato-LogStruct}, and section 8.1 of \cite{Kato-ToricSing}. However this converse is also implied by Proposition \ref{prop2} below, which is a stronger assertion. Proposition \ref{prop2} can be proved, adapting the argument in section 3.13 of \cite{Kato-LogStruct}. Because this argument is phrased in the framework of logarithmic structures on schemes, we give below an elementary self-contained proof of Proposition \ref{prop2} which does not use logarithmic geometry.
\end{remark}

\begin{proposition}\label{prop2} Let $U_X \subset X$ and $U_B \subset B$ be strict toroidal embeddings of varieties over $k$, and $x$ a closed point of $X$. Let $f: X \rightarrow B$ be a dominant $k$-morphism mapping $U_X$ into $U_B$. Set $b := f(x)$. Let $(V_B,\varphi_B)$ be an \'etale chart for $U_B \subset B$ at $b$. Assume that $f$ is logarithmically smooth at $x$ with respect to $U_X \subset X$ and $U_B \subset B$. Assume also that $k(x)=k$ or that $k(x)$ is algebraically closed. \\
Then there exist
\begin{enumerate}
\item an \'etale k-morphism $\pi: X' \rightarrow X$,
\item a closed point $x'$ on $X'$ with $\pi(x')=x$ and $k(x')=k(x)$,
\item an \'etale chart $(V_{X'}, \varphi_{X'})$ for $U_{X'}:= \pi^{-1}(U_X) \subset X'$ at $x'$,
\item a toric morphism $g: V_{X'} \rightarrow V_B$,
\item a translation $t:V_B \rightarrow V_B$ by a $k$-rational point on the big torus of $V_B$,
\end{enumerate}
such that the following diagram of rational maps commutes
$$ \begin{array}{lclcl}
                    X' & \stackrel{\varphi_{X'}}{\longrightarrow} & V_{X'} \\
                    \downarrow \pi \\
                    X  & & \downarrow g\\
                    \downarrow f \\
                    B & \stackrel{t\circ\varphi_B}{\longrightarrow} &V_B
\end{array} $$
If $(V_X, \varphi_X)$ is any \'etale chart for $U_X \subset X$ at $x$ such that $\varphi_X(x)$ belongs to the closed orbit of $V_X$, then we can choose $V_{X'}=V_X$, with $\varphi_{X'}(x')$ in the orbit of $\varphi_X(x)$. And when $k(x)$ is algebraically closed we can moreover take for $t$ the identity.
\end{proposition}

\noindent \emph{Proof.} The proof consists of several steps.

\emph{Some reductions.} Let $(V_X, \varphi_X)$ be an \'etale chart for $U_X \subset X$ at $x$ such that $\varphi_X(x)$ belongs to the closed orbit of $V_X$. Replacing $B$ and $X$ by suitable open subvarieties we may suppose that $\varphi_B$ and $\varphi_X$ are defined and \'etale everywhere. We can assume that $\dim B = \dim X$ by replacing $f: X \rightarrow B$ by $f\times (h \circ \varphi_X): X \rightarrow B \times \A^{\dim X - \dim B}$, with
$h:V_X \rightarrow \A^{\dim X - \dim B}$ a general enough toric morphism so that $f\times (h \circ \varphi_X)$ is still dominant and logarithmically smooth at $x$. \\
Choose a basis $c_1, \dots, c_n$ for the $\Z$-module of characters on the big torus $T_B$ of $V_B$. We will denote $\varphi_B^*(c_i)$ again by $c_i$, for $i=1, \dots, n$.

\emph{Choosing character bases.} Choose a basis $z_1, \dots, z_r, z_{r+1}, \dots, z_m$ for the $\Z$-module of characters on the big torus $T_X$ of $V_X$ such that $z_{r+1}, \dots, z_m$ form a basis for the $\Z$-module of characters on $T_X$ that are defined and not vanishing at $\varphi_X(x)$. We will denote $\varphi_X^*(z_i)$ again by $z_i$, for $i=1, \dots, m$. Because $\dim B = \dim X$ we have $m = n$. \\
From Lemma \ref{lemma1} below it follows that for $j=1, \dots, n=m$ we can write
\[ f^*(c_j) = u_j z_1^{e_{j,1}}z_2^{e_{j,2}} \cdots z_n^{e_{j,n}}, \tag{1} \]
with the $u_j$ suitable units in $\cO_{X,x}$. Moreover we can choose $e_{j,r+1}, \dots, e_{j,n}$ arbitrarily if we adapt the $u_j$ to these choices, since $z_{r+1}, \dots, z_n$ are units in $\cO_{X,x}$.

\emph{Changing coordinates by Hensel's Lemma.} Let $J$ be the logarithmic jacobian matrix of $f$, i.e. the square matrix consisting of the coefficients expressing $f^*(\frac{dc_j}{c_j}), j=1, \dots, n,$ as $\cO_{X,x}$-linear combinations of $\frac{dz_i}{z_i}, i=1, \dots, n$.
We denote by $J(x)$ the square matrix over $k(x)$ obtained from $J$ by evaluation at $x \in X$. Because $f$ is logarithmically smooth at $x$ we have $\det J(x) \neq 0$. From (1) and Lemma \ref{lemma2} it follows that the first $r$ columns of $J(x)$ equal the first $r$ columns of the matrix $E:= (e_{j,i})_{j,i=1, \dots, n}$. Thus the first $r$ columns of $E$ are linearly independent. Since the last $n-r$ columns of $E$ can be chosen arbitrarily, we can choose these such that $\det(E)\neq 0$.\\
For $j=1, \dots, n$ we set $\lambda_j:=1$ if $k(x)$ is algebraically closed. Otherwise $k(x)=k$, by our assumption on $k(x)$, and then we set $\lambda_j:=u_j(x) \in k$. Hence for $j=1, \dots, n$ we can write
\[ f^*(c_j) = \lambda_j w_j z_1^{e_{j,1}}z_2^{e_{j,2}} \cdots z_n^{e_{j,n}}, \tag{2} \]
with the $w_j$ suitable units in $\cO_{X,x}$. Moreover $w_j(x)=1$ when $k$ is not algebraically closed. \\
Since $\det(E)\neq0$, it follows from Hensel's Lemma that there exist units $\epsilon_1, \dots, \epsilon_n$ in the henselization of $\cO_{X,x}$ such that for $j=1, \dots, n$
\[ w_j = \epsilon_1^{e_{j,1}}\epsilon_2^{e_{j,2}} \cdots \epsilon_n^{e_{j,n}}. \tag{3} \]
We will use the change of coordinates $z_i \leftarrowtail \epsilon_i z_i$, but to make this precise we have to go to an \'etale extension $X'$ of $X$ and use a new \'etale chart $(V_{X'}, \varphi_{X'})$.

\emph{Construction of $\pi: X' \rightarrow X$.} There exists an \'etale morphism $\pi: X' \rightarrow X$ and a closed point $x'$ on $X'$, with $\pi(x')=x$ and $k(x')=k(x)$, such that $\epsilon_1, \dots, \epsilon_n$ are units in $\cO_{X',x'}$. We may even assume that $X'$ is affine and that $\epsilon_1, \dots, \epsilon_n$ are units in the coordinate ring of $X'$. Note that $(V_X, \varphi_X \circ \pi)$ is an \'etale chart at $x'$ for the toroidal embedding $U_{X'}:= \pi^{-1}(U_X) \subset X'$, but we will need another chart $(V_{X'}, \varphi_{X'})$.

\emph{Construction of the chart $(V_{X'}, \varphi_{X'})$.} By multiplicativity, the assignment $z_i \mapsto \epsilon_i$ extends uniquely to a homomorphism $\epsilon$ from to group of characters of $T_X$ to $\Gamma(\cO_{X'},X')^\times$. Set $V_{X'}:=V_X$ and let $\varphi_{X'} : X' \rightarrow V_{X'}=V_X$ be the unique $k$-morphism with
\[\varphi_{X'}^*(z) = \epsilon(z) (\varphi_X \circ \pi)^*(z) ,  \tag{4} \]
for each character $z$ of $T_X$ (note that the characters of $T_X$ that are regular on $V_X$ generate the coordinate ring of $V_X$). Note that $\varphi_{X'}(x')$ belongs to the orbit of $\varphi_X(x)$ under the action of $T_X$. We show below that the pair $(V_{X'}, \varphi_{X'})$ is an \'etale chart at $x'$ for $U_{X'}\subset X'$.

\emph{Construction of translation $\, t$ and toric morphism $g$.} Let $t: V_B \rightarrow V_B$ be the translation by the $k$-rational point of $T_B$ on which the  character $c_j$ takes the value $\lambda_j^{-1}$, for $j=1, \dots, n$. \\
Finally, let $g: V_{X'} \rightarrow V_B$ be the toric rational map defined by
\[g^*(c_j) = z_1^{e_{j,1}}z_2^{e_{j,2}} \cdots z_n^{e_{j,n}},  \tag{5} \]
for $j=1, \dots, n$. We show below that $g$ is regular at each point of $V_{X'}$, i.e. $g$ is a morphism. \\
From (4), (3), and (2) it follows that
\[ (t \circ \varphi_B) \circ f \circ \pi = g \circ \varphi_{X'}.  \tag{6} \]
Thus the diagram in \ref{prop2} is indeed commutative.

\emph{The rational map $g$ is regular on $V_{X'}$.} To prove this it suffices to show that $g^*(c)$ is regular on $V_{X'}$, for each character $c$ of $T_B$ that is regular on $V_B$. From (6) it follows that $(g \circ \varphi_{X'})^*(c)$ is regular on $X'$, hence by (4) also $(\varphi_X \circ \pi)^*(g^*(c))$ is regular on $X'$. Thus $g^*(c) \in \cO_{V_{X'},\varphi_X(x)}$ because the homomorphism
$\cO_{V_{X'},\varphi_X(x)} \rightarrow \cO_{X',x'}$
induced by $\varphi_X \circ \pi : X' \rightarrow V_{X'} = V_X$ is faithfully flat, since
$\varphi_X \circ \pi$ is \'etale. Moreover $g^*(c)$ is a character of $T_X$, hence its divisor on $V_X$ is supported on $V_X \setminus T_X$.
Because $\varphi_X(x)$ belongs to the closed orbit of $V_X$, all irreducible components of $V_X \setminus T_X$ contain $\varphi_X(x)$. Since we know already that $g^*(c)$ is regular at $\varphi_X(x)$, we conclude that $g^*(c)$ is regular at each point of $V_{X'}$.

\emph{The pair $(V_{X'}, \varphi_{X'})$ is an \'etale chart at $x'$ for $U_{X'}\subset X'$.} For this it suffices to show that $\varphi_{X'}: X' \rightarrow V_{X'} = V_X$ is \'etale at $x'$, because $\epsilon(z)$ in formula (4) is a unit in $\Gamma(\cO_{X'},X')$. Since $f$ is logarithmically smooth at $x$, formula (6) implies that $g \circ \varphi_{X'}$ is logarithmically smooth at $x'$. Since $\dim V_B = \dim V_{X'}$, this implies (by the definition of logarithmically smooth) that $\varphi_{X'}$ is logarithmically smooth at $x'$ with respect to $U_{X'}\subset X'$ and $T_X \subset V_X$. Hence Lemma \ref{lemma3} below (with $X, \rho, \psi$ replaced by $X', \varphi_X \circ \pi, \varphi_{X'})$ implies that $\varphi_{X'}$ is \'etale at $x'$. \\
This terminates the proof of Proposition \ref{prop2}.  $\square$

\begin{remark}\label{remark2-counterExample}\rm Note that the assumption on $k(x)$ in the statement of Proposition \ref{prop2} is always satisfied if $k=\R$. The following counterexample shows that we cannot omit this assumption in Proposition \ref{prop2}. Let $X = \Spec(\Q[x,y,y^{-1}]/(y^2-x+1))$,
$U_X = X \setminus V(x)$, where $V(x)$ denotes the locus of $x=0$, $B = \Spec(\Q[z])$, $U_B = B \setminus V(z)$, and let $f: X \rightarrow B$ be given by $z= yx^4$. Let $b = (0) \in B$ and $a$ the unique point in $X$ with $f(a) = b$. The morphism $f$ is logarithmically smooth at $a$ with respect to $U_X \subset X$ and $U_B \subset B$. However the conclusion in Proposition \ref{prop2} (with $x$ replaced by $a$) does not hold. Indeed, there does not  exist a unit $u$ in the henselization of $\cO_{B,b}$, and a unit $v$ in the henselization of $\cO_{X,a}$, such that $zu = (xv)^4$. Otherwise $yu = v^4$, and taking values at $a$ we see that then $\sqrt{-1}$ could be written as $\alpha \beta^4$, with $\alpha \in \Q$ and $\beta \in \Q(\sqrt{-1})$. However, this is impossible.
\end{remark}

\begin{lemma} \label{lemma1} Let $U_X \subset X$ be a strict toroidal embedding of varieties over $k$ and let $x$ be a closed point of $X$. Let $(V,\varphi)$ be an \'etale chart for $U_X \subset X$ at $x$. Let $y \in \cO_{X,x}$ and assume that the divisor of $y$ is supported on $X \setminus U_X$ in some Zariski neighborhood of $x$ in $X$. Then there exists a character $c$ of the big torus $T$ of $V$ such that $y/\varphi^*(c)$ is a unit in $\cO_{X,x}$.
\end{lemma}

\noindent \emph{Proof.} This is very well known. From Lemma \ref{lemma0}, with $X, D$ replaced by $V, V \setminus T$, it follows that there exists a unit $z$ in $\cO_{V,\varphi(x)}$ such that $y/\varphi^*(z)$ is a unit in $\cO_{X,x}$. The ideal generated by $z$ in $\cO_{V,\varphi(x)}$ is invariant under the action of $T$, hence it is generated by characters of $T$. Since it is principal, Nakayama's Lemma yields that it is generated by one of these characters (cf. section 3.3 of \cite{Fulton}).  $\square$

\begin{lemma} \label{lemma2} Let $U_X \subset X$ be a strict toroidal embedding of varieties over $k$, and $x$ a closed point of $X$. Let $(V,\varphi)$ be an \'etale chart of $U_X \subset X$ at $x$, and set $v := \varphi(x)$. Denote by $T$ the big torus of $V$, and set $D_V := V \setminus T, D_X := X \setminus U_X$. Choose a basis $z_1, \dots, z_r, z_{r+1}, \dots, z_m$ for the $\Z$-module of characters on $T$, such that $z_{r+1}, \dots, z_m$ form a basis for the $\Z$-module of characters on $T$ that are defined and not vanishing at $v$. Let $y \in \cO_{X,x}$. Consider $dy$ and $\frac{dz_i}{z_i}$, for $i = 1, \dots, m$, as elements of $(\Omega_X^1(\log D_X))_x$, and write
\[ dy = \sum_{i=1}^m \, a_i \, \frac{dz_i}{z_i} \, , \]
with $a_i \in \cO_{X,x}$, for $i = 1, \dots, m$. Then $a_i(x)=0$, for $i = 1, \dots, r$.
\end{lemma}

\noindent \emph{Proof.} Since $\varphi$ is \'etale at $x$, the $\cO_{X,x}$-module $\Omega_{X,x}^1$ is generated by the pullbacks under $\varphi$ of the elements of $\Omega_{V,v}^1$. Thus we may suppose that $X=V$, $\varphi = $id, and $x=v$. We may also assume that $y$ is a character of $T$ that is regular on $V$, since these generate the coordinate ring of $V$. Moreover, replacing $V$ by a suitable open toric subvariety (on which $z_{r+1}, \dots, z_m$ are regular), we may also suppose that $v$ belongs to the closed orbit of $V$. \\
Let $A$ be the set of characters of $T$ that belong to the group generated by $z_1, \dots, z_r$ and are regular on $V$. Let $B$ be the group generated by $z_{r+1}, \dots, z_m$. Then $V_0 := \Spec k[A]$ is toric and $T_0 := \Spec k[B] \cong (\G_m)^{m-r}$. Moreover, since $v$ belongs to the closed orbit of $V$, we have an isomorphism $V \cong V_0 \times T_0$ induced by $a \otimes b \mapsto ab$ for $a \in A, b \in B$, and the closed orbit of $V_0$ consists of only one point. \\
Because the lemma is trivial when $V$ is a torus, we may assume that $r=m$, and hence that the closed orbit of $V$ consists of only one point, namely $v$. This implies that each nontrivial character of $T$, that is regular on $V$, vanishes at $v$. Thus $y(v)=0$, because we may assume that $y$ is such a character. Writing $y = z_1^{e_1}z_2^{e_2} \cdots z_m^{e_m}$, with $e_1, \dots, e_m \in \Z$, we have
$\frac{dy}{y} = \sum_{i=1}^m e_i \frac{dz_i}{z_i}$. This terminates the proof of the lemma. $\square$

\begin{remark}\label{remark3}\rm Assume the notation of Lemma \ref{lemma2} and denote by $C$ the orbit of $v$ under the action of $T$.
Then in $\Omega_{\varphi^{-1}(C),x}^1$ we have the equality $dy = \sum_{i=r+1}^m a_i \frac{dz_i}{z_i}$. This follows e.g. from the argument in the proof of Lemma \ref{lemma2}.
\end{remark}

The following lemma is (modulo the terminology of logarithmic geometry) a special case of Proposition 3.8 in \cite{Kato-LogStruct}.
\begin{lemma} \label{lemma3} Let $V$ be an affine toric variety over $k$ and denote its big torus by $T$. Let $X$ be an algebraic variety over $k$, $x$ a closed point of $X$, and $\rho: X \rightarrow V$ an \'etale k-morphism. In particular, $(V, \rho)$ is an \'etale chart at $x$ of the toroidal embedding $U_X := \rho^{-1}(T) \subset X$. \\
Let $\psi: X \rightarrow V$ be a dominant $k$-morphism from $X$ to $V$, mapping $U_X$ into $T$. Assume the following two conditions.
\begin{enumerate}
\item For each character $c$ of $T$,we have that $\rho^*(c)/\psi^*(c)$ is a unit in $\cO_{X,x}$.
\item The morphism $\psi$ is logarithmically smooth at $x$ with respect to $U_X \subset X$ and $T \subset V$.
\end{enumerate}
Then $\psi$ is \'etale at $x$.
\end{lemma}

\noindent \emph{Proof.} Instead of relying on Proposition 3.8 in \cite{Kato-LogStruct}, we give a self-contained proof that does not use logarithmic geometry. \\
Condition (1) implies that $\psi(x)$ belongs to the $T$-orbit of $\rho(x)$ in $V$. Hence, replacing $V$ by a suitable open toric subvariety and $X$ by a suitable open subvariety, we may assume that $\rho(x)$ and $\psi(x)$ belong to the closed orbit of $V$. We denote the closed orbit of $V$ by $C$ and set $v:=\rho(x) \in C, w:=\psi(x) \in C$. \\
To show that $\psi$ is \'etale at $x$, it suffices to prove that $\psi$ is unramified at $x$, because $\psi$ is dominant with integral source and normal target (see e.g. Theorem 3.20 in chapter 1 of \cite{Milne-EtaleCoh}. Hence it suffices to prove the following two claims.

\begin{itemize}
\item \emph{Claim 1.} The ideal of $\psi^{-1}(C)$ in $\cO_{X,x}$ is generated by elements $\psi^*(c)$ with $c$ in the ideal of $C$ in $\cO_{V,w}$.

\item \emph{Claim 2.} The morphism $\psi|_{\psi^{-1}(C)}: \psi^{-1}(C) \rightarrow V$ is unramified at $x$.
\end{itemize}

\emph{Proof of Claim 1.} The ideal of $C$ in the coordinate ring of $V$ is generated by the characters of $T$ that are defined and vanishing at $v$ (or equivalently at $w$, because $w$ belongs to the orbit of $v$). Hence condition (1) implies that $\psi^{-1}(C) = \rho^{-1}(C)$ locally at $x$. Thus, again by condition (1), in order to prove Claim 1, it suffices to show that the ideal of $\rho^{-1}(C)$ in $\cO_{X,x}$ is generated by the elements $\rho^*(c)$ with $c$ running over all characters of $T$ that are defined and vanishing at $v$. But this follows directly from the fact that $\rho$ is \'etale at $x$, because these characters generate the ideal $I$ of $C$ in $\cO_{V,v}$, and because $\cO_{V,v}/I$ is normal and hence $\rho^*(I)$ is prime.

\emph{Proof of Claim 2.} Clearly it suffices to show that the morphism $\psi^{-1}(C) \rightarrow C$ induced by $\psi$ is \'etale at $x$. We will show this by using the jacobian criterium. Note that $C$ is smooth, and that $\psi^{-1}(C)$ is smooth at $x$, because $\psi^{-1}(C) = \rho^{-1}(C)$ locally at $x$ (as we saw in the proof of Claim (1)). Choose a basis $z_1, \dots, z_r, z_{r+1}, \dots, z_m$ for the $\Z$-module of characters on $T$, such that $z_{r+1}, \dots, z_m$ form a basis for the $\Z$-module of characters on $T$ that are defined and not vanishing at $v$ (or equivalently at $w$). Note that $z_{r+1}, \dots, z_m$ are uniformizing parameters for $C$ at $v$, and also at $w$. Put
$$ x_i = \rho^*(z_i), \, \, x'_i = \psi^*(z_i), $$
for $i=1, \dots, m$. Note that $x_{r+1}, \dots, x_m$ are uniformizing parameters for $\psi^{-1}(C)$ at $x$, because $\psi^{-1}(C) = \rho^{-1}(C)$ locally at $x$, and because $\rho$ is \'etale. Thus, by the jacobian criterium, we have to prove that
\[\det \left(\frac{\partial x_i'}{\partial x_j}(x)\right)_{i,j=\,r+1, \dots, m} \neq 0. \tag{3} \]
From condition (1) it follows that there are units $\epsilon_i$ in $\cO_{X,x}$ such that for $i = 1, \dots, m$
\[x_i' = \epsilon_i x_i. \tag{4} \]
Hence
\[ \psi^*(\frac{dz_i}{z_i}) = \frac{d(\epsilon_i x_i)}{\epsilon_i x_i} = \frac{d\epsilon_i}{\epsilon_i} + \frac{dx_i}{x_i} \,.  \tag{5} \]
Let $J$ be the logarithmic jacobian matrix of $\psi$, i.e. the square matrix of the coefficients expressing $\psi^*(\frac{dz_i}{z_i})$, $i=1, \dots, m$, as $\cO_{X,x}$-linear combinations of $\frac{dx_j}{x_j}\,$, $j=1, \dots, m$. We denote by $J(x)$ the square matrix over $k(x)$ obtained from $J$ by evaluation at $x \in X$. Because of condition (2), we have that $\det J(x) \neq 0$. Applying Lemma \ref{lemma2} (with $y$ replaced by $\epsilon_i$) and using the fact that the $\epsilon_i$ are units, together with (5), we see that the matrix formed by the last $m-r$ rows and the first $r$ columns of $J(x)$ is zero. Hence the submatrix $J_0$ of $J$, formed by the last $m-r$ rows and the last $m-r$ columns, satisfies $\det(J_0)(x) \neq 0$. Note that for $i=r+1, \dots, m$, the $x_i$, and hence also the $x_i'$, are units in $\cO_{X,x}$.
Hence $\left(\frac{\partial x_i'}{\partial x_j}(x)\right)_{i,j=\,r+1, \dots, m}$ can be obtained from $J_0(x)$ by multiplying the $i$-th row of $J_0(x)$ by $x_i'(x) \neq 0$ and dividing the $j$-th column of $J_0(x)$ by $x_j(x) \neq 0$ (see Remark \ref{remark3}). This yields (3) and terminates the proof of the lemma. $\square$

\bibliographystyle{plain}
\bibliography{JD}

\end{document}